\documentclass[a4paper,12pt]{article}
\usepackage{amssymb}
\usepackage{amsmath}
\usepackage{amsthm}
\usepackage{graphicx}
\numberwithin{equation}{section}

\newcommand{\nuij}{\nu_{i\to j}}
\newcommand{\nuji}{\nu_{j\to i}}
\newcommand{\ep}{\varepsilon}

\newcommand{\va}{\varphi}
\newcommand{\ppp}{\partial}

\newcommand{\ddd}{\mbox{div}\thinspace}

\newcommand{\weight}{e^{2s\va}}

\newcommand{\sumN}{\sum_{i=1}^N}

\newcommand{\R}{\mathbb{R}}

\newcommand{\N}{\mathbb{N}}
\newcommand{\Z}{\mathbb{Z}}

\newcommand{\www}{\widetilde}

\newcommand{\ooo}{\overline}
\newcommand{\OOO}{\Omega}
\newcommand{\LLLL}{L^2(Q)}

\allowdisplaybreaks

\title{Carleman estimate with piecewise weight and applications to
inverse problems for first-order transport equations}

\author{
$^1\;$Piermarco Cannarsa, 
$^2\;$Giuseppe Floridia, 
$^3\;$Masahiro Yamamoto }

\date{}
\begin{document}
\maketitle

\begin{abstract}
We consider a first-order transport equation 
$\ppp_tu(x,t) + (H(x)\cdot\nabla u(x,t)) + p(x)u(x,t) = F(x,t)$
for $x \in \OOO \subset \R^d$, where $\OOO$ is a bounded domain and $0<t<T$.
We prove a Carleman estimate for more generous condition on the 
principal coefficients $H(x)$ than in the existing works.
The key is the construction of a piecewise smooth weight 
function in $x$ according to a suitable decomposition of $\OOO$.   
Our assumptions on $H$ generalize the conditions in the existing articles,
and require that a directed graph created by the corresponding 
stream field has no closed loops.  Then, we apply our Carleman estimate 
to two inverse problems of 
determinination of an initial value and one of a spatial factor of 
a source term, so that we establish Lipschitz stability estimates for the 
inverse problems.
\end{abstract}

\section{Introduction and motivation}

We consider a first-order transport equation:
$$
\ppp_tu(x,t) + (H(x)\cdot\nabla u(x,t)) + p(x)u(x,t) = F(x,t), \quad 
0<t<T, \, x\in \OOO,
                                   \eqno{(1.1)}
$$
where $\OOO \subset \R^d$ is a bounded domain with smooth boundary 
$\ppp\OOO$.    
Throughout this article, we assume
$$
H = (h_1, ..., h_d) \in C^1(\ooo{\OOO}), \quad
\vert H(x)\vert \ge \delta_0 > 0, \quad x\in \ooo{\OOO},        \eqno{(1.2)}
$$
where $\delta_0 > 0$ is a constant, and let 
$p\in L^{\infty}(\OOO)$.  Here and henceforth,
let $x:= (x_1, ..., x_d) \in \R^d$, $\ppp_k := \frac{\ppp}{\ppp x_k}$ for 
$k=1, ..., d$ and $\nabla := (\ppp_1, ..., \ppp_d)^T$.  Moreover,
$(v\cdot w)$ is the scalar product of $v, w\in \R^d$, and $^T$ denotes the 
transpose of vectors and matrices under consideration.

By $\nu = \nu(x)$ we denote the outward unit normal vector 
to $\OOO$ at $x$, and set 
$$
\left\{ \begin{array}{rl}
&Q:= \OOO \times (0,T), \quad S^{d-1}:= \{ x\in \R^d;\,
\vert x\vert = 1\},               \cr\\
& \ppp\OOO_+:= \{ x\in \ppp\OOO;\, (H(x)\cdot \nu(x)) \ge 0\}, \cr\\
& \ppp\OOO_-:= \{ x\in \ppp\OOO;\, (H(x)\cdot \nu(x)) < 0\}. 
\end{array}\right.
                            \eqno{(1.3)}
$$

For (1.1), we consider
\\
{\bf Observability:}   
{\it 
Let $F=0$ in $Q$.  Determine $u(\cdot,0)$ in $\OOO$ by 
$u\vert_{\ppp\OOO\times (0,T)}$.}
\\
{\bf Inverse source problem:} 
{\it 
Let $F(x,t)= R(x,t)f(x)$ in (1.1) with a suitably given function $R(x,t)$ for 
$(x,t) \in Q$.  Then determine $f(x)$, $x \in \OOO$ by data
$(u\vert_{\ppp\OOO\times (0,T)}, \, \, u(\cdot,0)\vert_{\OOO})$.}
\\

The inverse source problem with $x$-dependent $R(x,t)$ is essential
for inverse coefficient problems such as the determination of the coefficient
$p(x)$ in $\ppp_tu(p) + (H\cdot \nabla u(p)) + p(x)u = 0$. Indeed, 
subtracting the equations for $p$ and $q$ and setting $f:= p-q$, we can 
reduce the inverse coefficient problem to an inverse source problem
(e.g., \cite{Y24}).

Equations of type (1.1) appear in mathematical physics, for example:  
$$
\mbox{conservation law:} \quad 
\ppp_tu(x,t) + \mbox{div}\, (V(x)u(x,t)) = 0 \quad \mbox{in $Q$}, 
                               \eqno{(1.4)}
$$
where $V:= (V_1,V_2,V_3)^T$.
\\
In governing equations for a linearized compressible non-viscous fluid
(called the Euler equations):
$$
\left\{ \begin{array}{rl}
& \ppp_tV(x,t) + (H(x,t)\cdot \nabla)V(x,t) + \nabla p(x,t) = 0,\\
& \mbox{div}\, V(x,t) = 0 \quad \mbox{in $Q$},
\end{array}\right.
\eqno{(1.5)}
$$
where $V = (V_1,V_2, V_3)^T$ is the velocity field and 
$p$ is the pressure, and 
we set $(H\cdot \nabla) V:= \left( \sum_{k=1}^3 h_k\ppp_k V_i\right)
_{1\le i \le 3}$.
\\

For discussing inverse problems, our first key is a weighted estimate 
called a Carleman estimate, and the second key is 
the methodology based on Carleman estimates by Bukhgeim and Klibanov 
\cite{BK}. For studying first-order partial differential equations, 
the method of characteristics is also a standard tool, 
but it is difficult to obtain $L^2$-estimates by 
characteritics, and moreover, is strictly limited 
to first-oder equations.  For example, such a technique is 
not applicable to the compressible fluid system where (1.5) is coupled 
with a parabolic equation of second order.
However, our method by Carleman estimates can directly work 
(e.g., Imanuvilov and Yamamoto \cite{IY20}). 

As for such an approach by Carleman estimates for first-order 
transport equations,
we refer, for example, to Cannarsa, Floridia, G\"olgeleyen, Yamamoto
\cite{CFGY}, Gaitan and Ouzzane \cite{GO}, G\"olgeleyen and Yamamoto \cite{GY},
Chapter 3 in Yamamoto \cite{Y24}.
In those works, some conditions are assumed on $H(x)$ in order to prove
Carleman estimates.  We can paraphrase their conditions as 
\\
{\bf Condition (A)}
$$
\left\{ \begin{array}{rl}
& \mbox{there exist a constant vector $v\in S^{d-1}$ and 
a constant $\delta_1 > 0$}\\
& (H(x)\cdot v) \ge \delta_1 \quad \mbox{for all $x\in \ooo{\OOO}$.}
\end{array}\right.
                                     \eqno{(1.6)}
$$

This condition means that the range of the solid angles of $H(x)$ over
$x \in \ooo{\OOO}$ is included in the half $\frac{\vert S^{d-1}\vert}{2}$
of the unit sphere.  Here, recalling that $\vert H(x)\vert > 0$ for 
$x\in \ooo{\OOO}$, we define the solid angle of $H(x)$ by 
$\frac{H(x)}{\vert H(x)\vert}$.

Under Condition (A), we can prove a Carleman estimate and directly
derive an observabilty inequality.
\\
{\bf Proposition 1.1 (Carleman estimate under Condition (A))}
\\
{\it
We choose 
$$
\left\{ \begin{array}{rl}
& \va(x,t):= \vert x+rv\vert^2 - \beta t, \quad x\in \OOO, \, 
0<t<T, \\
& r > \frac{1}{2\delta_1}(2(\max_{x\in \ooo{\OOO}}\vert x\vert)
\Vert H\Vert_{C(\ooo{\OOO})} + \beta).
\end{array}\right.
                                         \eqno{(1.7)}
$$
Then, there exist constants $s_0>0$ and $C>0$ such that 
$$
s\int_{\OOO} \vert u(x,0)\vert^2 e^{2s\va(x,0)} dx
+ s^2\int_Q \vert u\vert^2 e^{2s\va(x,t)} dxdt
+ s\int_{\ppp\OOO_-\times (0,T)} \vert u\vert^2 e^{2s\va(x,t)}dS dt
$$
$$
\le C\int_Q \vert \ppp_tu + (H\cdot\nabla u)\vert^2 e^{2s\va(x,t)} dxdt
+ Cs\int_{\ppp\OOO_+\times (0,T)} \vert u\vert^2 e^{2s\va(x,t)}dS dt
$$
$$
+ Cs\int_{\OOO} \vert u(x,T)\vert^2 e^{2s\va(x,T)} dx
                                      \eqno{(1.8)}
$$
for all $s\ge s_0$ and $u \in H^1(Q)$.}
\\
For completeness, we provide the proof of 
Proposition 1.1 in Section 6.

It is known that in general, Carleman estimates produce several kinds of 
stability estimates 
for inverse problems (e.g., Bellassoued and Yamamoto \cite{BY}, 
Bukhgeim and Klianov \cite{BK}, Klibanov \cite{Kl}, Klibanov and Timonov 
\cite{KT}, Yamamoto \cite{Ya}, \cite{Y24}).  We here show one application. 
\\ 
{\bf Proposition 1.2 (Observability inequality under Condition (A))}
\\
{\it
We assume that 
$$
T > \frac{2(\max_{x\in \ooo{\OOO}}\vert x\vert)(\delta_1 + \Vert H\Vert
_{C(\ooo{\OOO})} )} {\delta_1^2}.
$$
Then, we can find a constant $C>0$ such that 
$$
\Vert u(\cdot,0)\Vert_{L^2(\OOO)} \le C\Vert u\Vert
_{L^2(\ppp\OOO\times (0,T))}.
$$
}

In Section 4.1, we describe the argument how to 
derive the observability inequality from a relevant Carleman estimate
for a more general case than the case where Condition (A) holds, and 
so we can see that Proposition 1.1 yields Proposition 1.2 as a special case.

See Komornik \cite{Ko} for examples of observability inequalities for 
wave equations and the method called the multiplier method, which is different
from ours, is based on Carleman estimates.
 
Condition (A) is quite restrictive.  In other words, we can prove  
Carleman estimates in several cases where Condition (A) is not satisfied.
The main purpose of this article is to generalize Condition (A),  
prove a Carleman estimate for a wider class of $H(x)$, and apply 
such an estimate to inverse problems.

We can refer, for example, to Isakov's monograph \cite{Is} for
various inverse problems.  For similar inverse problems for 
first-order transport equations, see Cannarsa, Floridia and Yamamoto
\cite{CFY}, Floridia and Takase \cite{FT1, FT2}, 
Floridia, Takase and Yamamoto \cite{FTY}.
H\"olmander \cite{Ho} provides a general theory for Carleman esimates, but
we need piecewise weights, and his theory is not directly applicable
to inverse problems.  Therefore, our direct approach is more effective.
  
This article is composed of six sections.  In Section 2, we generalize 
Condition (A) and state the  main Carleman estimate.
In Section 3, we prove
it.  Section 4 is devoted to applications to observability 
and an inverse source problem under our relaxed conditions on $H$.  
In Section 5, we give concluding 
remarks and, in Section 6, we prove Proposition 1.1 and another 
Carleman estimate (Proposition 2.1) stated in Section 2.

\section{Main Carleman estimate} 

In order to derive a Carleman estimate, Condtion (A) is sufficient
but not a necessary condition.  However, if Condition (A) does not hold, 
then we may not be able to prove a Carleman estimate which is applicable to 
inverse problems or some uniqueness results related to a solution $u(x,t)$ to 
(1.1), as the following example shows.
\\
{\bf Example 2.1.}
\\
We consider the two-dimensional case: $d=2$.
Let $\OOO := \{ (x_1,x_2)\in \R^2;\, 1 < \vert x\vert < 2\}$ and
$H(x_1,x_2) := (-x_2,\, x_1)^T$.  We consider
$$
u_0(x_1,x_2,t) = (\vert x\vert^2-1)(\vert x\vert^2 - 4), \quad 
(x_1,x_2,t) \in Q:= \OOO \times (0,T).
$$
Then, 
$$
\left\{ \begin{array}{rl} 
& \ppp_tu_0(x,t) + (H(x)\cdot \nabla u_0(x,t)) = 0 \quad
\mbox{in $\OOO \times (0,\infty)$}, \\
& u_0 = 0 \quad \mbox{on $\ppp\OOO\times (0,\infty)$},
\end{array}\right.
$$
which means that the conclusion of Proposition 1.2 does
not hold.  Indeed,  
$$
\left\{ \frac{H(x)}{\vert H(x)\vert};\, x\in \ooo{\OOO}
\right\} 
= \left\{ \frac{(-x_2,x_1)}{\vert x\vert}; \, 1 < \vert x\vert < 2\right\}
= S^1,
$$
that is, (1.6) is not satisfied for any $v\in S^1$.
If we have a relevant Carleman estimate, then we can apply such an  
estimate to conclude that $ u_0=0$ in $Q$ similarly to 
Proposition 1.2.  Therefore, by such $u_0$, we see that we cannot
have any Carleman estimates in this case.

We emphasize that condition (1.6) is not a necessary condition in order to 
establish an applicable Carleman estimate to inverse problems.
Indeed, we consider a case where $H(x)$ is given by a scalar potential 
function:
$$
\OOO:= \{ (x_1,x_2) \in \R^2;\, 1 < \vert x\vert < 2\}, \quad 
H(x) := (2x_1, 2x_2) = \nabla (\vert x\vert^2).
$$
Clearly $2(x_1,x_2)$ does not satisfy (1.6) in $\OOO$.  However we can 
choose a weight function to establish a Carleman estimate.

We can state this case as
\\
{\bf Proposition 2.1 (Carleman estimate without Condition (A))}
{\it
We assume 
$$
H(x) = \nabla \rho(x), \quad x\in \ooo{\OOO},    \eqno{(2.1)}
$$
with some $\rho\in C^2(\ooo{\OOO})$ satisfying
$\vert \nabla \rho\vert > \beta> 0$ on $\ooo{\OOO}$.
We set $\va(x,t):= \rho(x) - \beta t$ for $x\in \OOO$ and 
$0<t<T$.  Then Carleman estimate (1.8) holds.
}

Condition (2.1) means that we consider a potential flow and 
can validate a Carleman estimate although (1.6) is not satisfied.
In this article, not assuming such special forms of $H$ such as (2.1),
we generalize (1.6) to establish a Carleman estimate.

For the statement of our main result, we introduce notations and definitions.
The main idea for relaxing conditions on $H$,
is to use a partition of $\OOO$ such that, in each subdomain, condition (1.6) is satisfied.  The proof of Carleman 
estimate over $\OOO$ can be reduced to estimations in each 
subdomain, but for synthesizing all the integrals in the subdomains, it is 
necessary to absorb all the integral terms on the interfaces
created by the partition.  Thus we need to assume some conditions,
which can be described as a directed graph.

We introduce a partition of $\OOO$:
$$
\ooo{\OOO} = \bigcup_{i=1}^N \ooo{\OOO_i},
$$
where $\OOO_i \subset \OOO$, $1\le i \le N$,
are non-empty disjoint subdomains.

By $\vert \gamma_{ij}\vert$, we denote 
the measure of the surface $\gamma_{ij}$.  
We assume 
$$
\mbox{$\gamma_{ij}:= \ppp\OOO_i \cap \ppp\OOO_j$ with
$i\ne j$ are of piecewise $C^2$-class if $\vert \gamma_{ij}\vert > 0$}.
                                   \eqno{(2.2)}
$$
By $\nuij$ we denote the unit normal vector to 
$\gamma_{ij}$ directed from 
$\OOO_i$ towards $\OOO_j$.  Note that $\nuij = -\nuji$.

We formulate 
\\
{\bf Condition (B).}
\\
{\it
The partition $\{ \OOO_i\}_{1\le i \le N}$ satisfies:
\\
(i) On each $\gamma_{ij}$, one of the following three cases hold:
$$
\left\{ \begin{array}{rl}
& (H\cdot \nuij) > 0 \quad \mbox{almost everywhere
on $\gamma_{ij}$},\\
\mbox{or} & (H\cdot \nuij) < 0 \quad \mbox{almost everywhere
on $\gamma_{ij}$},\\
\mbox{or} & (H\cdot \nuij) = 0 \quad \mbox{almost everywhere
on $\gamma_{ij}$}.
\end{array}\right.                          
$$
\\
(ii) For each $i=1,..., N$, there exist a constant $\delta_1 > 0$ and
a constant vector $v_i \in S^{d-1}$ such that 
$$
\min_{1\le i \le d} (H(x)\cdot v_i) > \delta_1 \quad 
\mbox{for all $x \in \ooo{\OOO_i}$}.
$$
}
\\
{\bf Remark.}  We do not consider the case where
$(H\cdot \nuij)$ changes signs infinitely many times, which 
makes point (i) of Condition (B) impossible.

Henceforth, to simplify the description, we neglect $\gamma_{ij}$ if 
$\vert \gamma_{ij} \vert = 0$. 

Next we correspond a directed graph $\Lambda$ to the partition 
$\{ \OOO_i\}_{1\le i\le N}$ as follows.
\\
(I) We correspond a node point $O_i$ to $\OOO_i$ for 
$1\le i \le N$.
\\
(II) By a segment $\Gamma_{ij}$, we connect two distinct nodes $O_i$ and
$O_j$ with $i\ne j$ if and only if $\vert \gamma_{ij} \vert 
= \vert \ppp\OOO_i\cap \ppp\OOO_j\vert > 0$ and 
$(H\cdot \nuij) \ne 0$ on $\gamma_{ij}$.  Moreover the segment 
$\Gamma_{ij}$ is directed from $O_i$ to $O_j$ if 
$(H\cdot \nuij) > 0$ almost everywhere on $\gamma_{ij}$, 
and from $O_j$ to $O_i$ if $(H\cdot \nuij) < 0$
almost everywhere on $\gamma_{ij}$.

We note that if $\vert \gamma_{ij} \vert = 0$ or
$(H\cdot \nuij) = 0$ on $\gamma_{ij}$, then 
$O_i$ and $O_j$ are not connected by a directed segment.

With the nodes $\{ O_i\}_{1\le i \le N}$ and the 
directed segments $\{ \Gamma_{ij}\}$, we compose a 
directed graph $\Lambda$.
Since $H(x)$ describes the stream field, we can interpret
$\Lambda$ as a chart of the stream field of (1.1). 
\\
{\bf Definition 2.1.}
\\
{\it 
We call that a directed graph $\Lambda$ has a closed loop if we can 
choose a subset $\{ \Gamma_{i_1j_1}\}$ of all the  
directed segments $\{ \Gamma_{ij}\}$ such that 
$\bigcup_{i_1,j_1} \Gamma_{i_1j_1}$ is a oriented closed Jordan curve. 
}

Then, we can formulate
\\
{\bf Condition (C).}
\\
{\it The directed graph $\Lambda$ associated with the partition $\{ \OOO_i\}
_{1\le i \le N}$ has no closed loops.
}

We set 
$$
Q_i:= \OOO_i \times (0,T), \quad 1\le i \le N.
$$
Now we are ready to state our main result.
\\
{\bf Theorem 2.1 (Carleman estimate)}
\\
{\it
We assume Conditions (B) and (C).  Then we can find constants 
$r_i > 0$ with $1\le i \le N$ and $\beta > 0$, and set 
$$
\va_i(x,t) = \vert x + r_iv_i\vert^2 - \beta t \quad 
\mbox{in $Q_i$},                                            \eqno{(2.3)}
$$
such that 
$$
s\sumN \int_{\OOO_i} \vert u(x,0)\vert^2 e^{2s\va_i(x,0)} dx
+ s^2\sumN \int_{Q_i} \vert u\vert^2 e^{2s\va_i} dxdt
+ s\int_{\ppp\OOO_- \times (0,T)} \vert u\vert^2 e^{2s\va} dS dt 
                                         \eqno{(2.4)}
$$
\begin{align*}
\le& C\sumN \int_{Q_i} \vert \ppp_tu + (H\cdot \nabla u) + pu\vert^2 
e^{2s\va_i} dxdt \\
+ & Cs \int_{\ppp\OOO_+\times (0,T)} \vert u\vert^2 e^{2s\va} dSdt 
+ Cs \sumN \int_{\OOO_i} \vert u(x,T)\vert^2 e^{2s\va_i(x,T)}dx
\end{align*}
for all $s > s_0$ and $u\in H^1(Q)$.
}
\\

Related to the existence of closed loops in directed graphs, we discuss 
the following \\
{\bf Example 2.2.}
Let $x = (x_1,x_2) \in \R^2$ and $\OOO:= \{x; \vert x\vert < 1\}$.
We use the two-dimensional polar coordinates:
$$
x_1 = r\cos \theta, \quad x_2 = r\sin \theta
$$
with $r\ge 0$ and $0 \le \theta < 2\pi$.
We consider a function $p(\theta)$ such that 
$$
p \in C^1[0,2\pi], \quad p(0) = 0, \quad 
p(2\pi) = 2\pi m, \quad p'(2\pi) = p'(0)                         \eqno{(2.5)}
$$
with $m \in \N \cup \{ 0\}$.
and  
$$
H(x_1,x_2) = (-\sin p(\theta), \, \cos p(\theta)),\quad 
0 \le \theta < 2\pi.
$$
Then we see that $H \in C^1(\ooo\OOO)$.  If the range $\{ p(\theta);\,
0 \le \theta < 2\pi\}$ is wider than $(0,\pi)$, then (1.6) fails.
However, we may be able to establish a Carleman estimate in many cases.
We consider special partitions: 
$$
\OOO_i:=\{ \theta_{i-1} < \theta < \theta_i \},\quad 1\le i \le N
$$
where $0=: \theta_0 < \theta_1 < \cdots < \theta_{N-1} < \theta_N:= 2\pi$
and $\OOO_{N+1} := \OOO_1$.  
Then we can correspond the nodes as $O_i := (\cos \theta_i,\, \sin \theta_i)$ 
for $0\le i \le N$.  Then,  
$$
\gamma_{i,i+1} := \ppp\OOO_i \cap \ppp\OOO_{i+1}
= \{\theta=\theta_i \}, \quad 1\le i \le N.
$$

Henceforth we set 
$$
q(\theta): = p(\theta) - \theta, \quad 0\le \theta < 2\pi.
$$
Then,  
$$
\left\{ \begin{array}{rl}
& \nu_{i \to i+1} = (-\sin \theta_i, \, \cos \theta_i)^T, \\
& (H \cdot \, \nu_{i \to i+1})  
= \sin p(\theta_i) \sin \theta_i + \cos p(\theta_i) \cos \theta_i
= \cos q(\theta_i) \quad \mbox{for $1\le i\le N$}.
\end{array}\right.
                                                \eqno{(2.6)}
$$

Assigning a node $O_i$ to $\OOO_i$ for $1 \le i \le N$ and setting
$O_{N+1} = O_1$, we connect $O_i$ to $O_{i+1}$ with the direction 
$O_i \rightarrow O_{i+1}$ if
$\cos q(\theta_i) > 0$, and $O_{i+1} \rightarrow O_i$ if 
$\cos q(\theta_i) < 0$, and do not connect such nodes 
if $\cos q(\theta_i) = 0$.  Hence,  
we construct a directed graph $\Lambda$ with nodes $O_1, ..., O_N$ and 
$N$-directed segments $\Gamma_{12}, \Gamma_{23}, ...., \Gamma_{N-1,N},
\Gamma_{N, N+1}$.

Then, we can prove
\\
{\bf Lemma 2.1.}
\\
{\it
{\bf (i) Case $m=0$ or $m\ge 2$ in (2.5)}.
\\
We can find a partition $0<\theta_1 < \cdots < \theta_{N-1}
< 2\pi$ such that 
$$
\vert \{ p(\theta);\, \theta_{i-1} < \theta < \theta_i \}\vert < \pi 
\quad \mbox{for $1\le i \le N$}            \eqno{(2.7)}
$$
and the corresponding directed graph $\Lambda$ has no closed loops.
\\
{\bf (ii) Case $m=1$ in (2.5)}.
\\
If we can find $\ell \in \Z$ such that 
$$
\mbox{$\frac{\pi}{2} + \ell\pi$ is an interior point in
$\{ q(\theta);\, 0 < \theta < 2\pi\}$},        \eqno{(2.8)}
$$
then the conclusion of case (i) holds.
}

We consider $p(\theta) = m\theta$ with $m\in \N \cup \{ 0\}$ as
a function satisfying (2.5).  Lemma 2.1 implies that:
If $m=0$ or $m\ge 2$, then we can find a directed graph without closed 
loops, and if $m=1$, then any directed graph has a closed loop.
Example 2.1 indicates a rare case where the graph has a closed loop.
Moreover, in the case $m=1$, we recall that 
we cannot obtain any Carleman estimate
which is applicable to the uniqueness of the solution $u$ by 
boundary measurements. 

Moreover, we note that if $H$ rotates fast in the sense that $m\ge 2$, 
then we can construct a directed graph without closed loops.
We note that, for $m\ge 2$, condition (1.6) does not hold.
Thus Theorem 2.1 yields a Carleman estimate under Conditions 
(B) and (C) which are weaker than Condition (A).
\\
\vspace{0.2cm}
\\
{\bf Proof of Lemma 2.1.}
\\
{\bf Case (ii).}
By (2.8), we can choose $\ell \in \Z$, $\theta_*, \theta^* \in 
(0, 2\pi)$ and sufficiently small $\ep>0$ such that 
$$
\frac{\pi}{2} + \ell \pi - \ep < q(\theta_*) < \frac{\pi}{2} + \ell \pi 
< q(\theta^*) < \frac{\pi}{2} + \ell \pi + \ep.                \eqno{(2.9)}
$$
Choosing $0=: \theta_0 < \theta_1 < \cdots < \theta_{N-1} < \theta_N:= 2\pi$ 
such that 
$$
\theta_{n_0}:= \theta_* < \theta^* := \theta_{n_0+1}
$$
with some $n_0 \in \{0, 1, ..., N-1\}$, by (2.9) we have $(\cos q(\theta_*))
(\cos q(\theta_*)) < 0$.  Hence (2.6) implies that 
$(H\cdot \nu_{n_0\to n_0+1})$ and $(H\cdot \nu_{n_0+1\to n_0+2})$ have 
different signs.
Moreover, we choose $\theta_i$ such that $\max_{0\le i \le N-1}(\theta_{i+1}
- \theta_i)$ is sufficiently small.  Then the continuity of $p$ implies 
(2.7).  
In other words, $\Gamma_{n_0, n_0+1}$ and $\Gamma_{n_0+1,n_0+2}$ are 
oppositely directed. Thus the proof of Case (ii) is finished.
\\
{\bf Case (i).}
The intermediate value theorem implies
$$
\{ q(\theta);\, 0\le \theta \le 2\pi\} \supset [q(0), \, q(2\pi)]
\quad \mbox{if $q(2\pi)\ge 0$}
$$
and
$$
\{ q(\theta);\, 0\le \theta \le 2\pi\} \supset [q(2\pi), \, q(0)]
\quad \mbox{if $q(2\pi)\le 0$}.
$$
Therefore,
$$
\{ q(\theta);\, 0\le \theta \le 2\pi\} \supset 
\left\{ \begin{array}{rl}
& [-2\pi,\, 0] \quad \mbox{if $m=0$}, \\
& [0, \, 2(m-1)\pi] \supset [0, \, 2\pi] 
\quad \mbox{if $m \ge 2$}.
\end{array}\right.
$$
Hence,
$$
\{ q(\theta);\, 0\le \theta \le 2\pi\} \ni 
\left\{ \begin{array}{rl}
& -\frac{\pi}{2} \quad \mbox{if $m=0$}, \\
&  \frac{\pi}{2} \quad \mbox{if $m \ge 2$}.
\end{array}\right.
$$
Since (2.8) is satisfied, similarly to Case (ii), we can verify the 
lemma in Case (i).  Thus the proof of Lemma 2.1 is complete.
$\blacksquare$
\section{Proof of Theorem 2.1}

Henceforth, for brevity, we set
$Q_i:= \OOO_i \times (0,T)$, and, for given $\va_i \in C^1(\ooo{Q_i})$ 
with $1\le i \le N$, we define a piecewise $C^1$-function on $\ooo{Q}
= \ooo{\OOO} \times [0,T]$ by 
$$
\va(x,t) = \va_i(x,t) \quad \mbox{for $(x,t) \in Q_i$ and 
$1\le i \le N$}.
$$
Also, for $U, V\in L^1(Q)$, we write
\begin{align*}
& \int_Q Ue^{2s\va} dxdt = \sum_{i=1}^N \int_{Q_i} Ue^{2s\va_i} dxdt,\\
& \int_{\OOO} Ve^{2s\va(x,\tau)} dx = \sum_{i=1}^N \int_{\OOO_i} 
Ve^{2s\va_i(x,\tau)} dx,\quad \mbox{$\tau = 0$ or $=T$}, 
\end{align*}
if there is no fear of confusion.  Moreover let 
$\nu_{\OOO_i}$ denote the outward normal vector to $\ppp\OOO_i$.

Henceforth, we recall that $C>0$ denotes a generic constant which is 
independent of
$s$ but may dependent on other quantities such as $\OOO, T, H, p$.

It is sufficient to prove Theorem 2.1 by assuming $p\equiv 0$ in (1.1).
Indeed, let Theorem 2.1 be proved for case $p\equiv 0$ in $Q$.
Then,
\begin{align*}
& s\int_{\OOO} \vert u(x,0)\vert^2 e^{2s\va(x,0)} dx
+ s^2 \int_{Q} \vert u\vert^2 e^{2s\va} dxdt
+ s\int_{\ppp\OOO_- \times (0,T)} \vert u\vert^2 e^{2s\va} dS dt \\
\le& C \int_Q \vert F - pu\vert^2 e^{2s\va} dxdt 
+ Cs\int_{\ppp\OOO_+ \times (0,T)} \vert u\vert^2 e^{2s\va} dSdt 
+ Cs \int_{\OOO} \vert u(x,T)\vert^2 e^{2s\va(x,T)}dx
\end{align*}
for all $s > s_0$ and $u\in H^1(Q)$.  Therefore,
\begin{align*}
& s\int_{\OOO} \vert u(x,0)\vert^2 e^{2s\va(x,0)} dx
+ (s^2 - C\Vert p\Vert_{L^{\infty}(\OOO)}^2)
\int_{Q} \vert u\vert^2 e^{2s\va} dxdt\\
+ &s\int_{\ppp\OOO_- \times (0,T)} \vert u\vert^2 e^{2s\va} dS dt \\
\le& C \int_Q \vert F \vert^2 e^{2s\va} dxdt 
+ Cs\int_{\ppp\OOO_+ \times (0,T)} \vert u\vert^2 e^{2s\va} dSdt\\ 
+ & Cs \int_{\OOO} \vert u(x,T)\vert^2 e^{2s\va(x,T)}dx
\end{align*}
for all $s > s_0$.

Hence, choosing $s_0 > 0$ sufficiently large such that 
$s_0 > \sqrt{2C}\Vert p\Vert_{L^{\infty}(\OOO)}$ for 
example, we can replace the second term on the left-hand side
by 
$$
C_1s^2 \int_{Q} \vert u\vert^2 e^{2s\va} dxdt.
$$
Thus the proof for $p\equiv 0$ is sufficient.
\\

We set $w_i:= e^{s\va_i}u_i$,
$$
P_iw_i:=e^{s\va_i}(\ppp_t + (H\cdot \nabla))(e^{-s\va_i}w)
$$
for $1\le i \le N$.  We recall that $H(x) = (h_1(x), ..., h_d(x))$ and 
we write $w := w_i$ and $Pw := P_iw_i$ in $Q_i:= \OOO_i \times (0,T)$ for 
$l \le i \le N$.  Let $\Vert \cdot\Vert$ stand for $\Vert \cdot\Vert_{L^2(Q)}$. We set $(U, \, V) = \int_Q U(x,t)V(x,t) dxdt$, if not specified extra.

Then, direct calculations yield
\begin{align*}
&Pw(x,t) := e^{s\va(x,t)} (\ppp_t + (H\cdot \nabla))(e^{-s\va(x,t)}w(x,t))\\
= & \ppp_tw(x,t) + (H(x)\cdot \nabla w(x,t)) 
- s(\ppp_t\va + (H\cdot \nabla\va))w(x,t)
\end{align*}
$$
=: \ppp_tw(x,t) + (H(x)\cdot \nabla w(x,t)) - sB(x,t)w(x,t) \quad
\mbox{in $Q$},
$$
where 
$$
\left\{ \begin{array}{rl}
& B_i(x,t):= \ppp_t\va_i(x,t) + (H(x)\cdot \nabla \va_i(x,t)), \quad (x,t) 
\in Q_i, \\
& B(x,t) = B_i(x,t) \quad \mbox{in $Q_i$ for $1\le i \le N$}.
\end{array}\right.
                                                  \eqno{(3.1)}
$$
Then,
\begin{align*}
& \int_Q \vert Pw\vert^2 dxdt = \Vert Pw\Vert^2_{L^2(Q)}
= \Vert \ppp_tw + (H\cdot\nabla w) - sBw\Vert^2_{L^2(Q)}\\
=& \Vert \ppp_tw + (H\cdot\nabla w)\Vert^2_{L^2(Q)}
+ s^2\Vert Bw\Vert^2_{L^2(Q)}
- 2s\int_Q (\ppp_tw + (H\cdot\nabla w))Bw dxdt
\end{align*}
$$
\ge s^2\Vert Bw\Vert^2_{L^2(Q)}
- 2s\int_Q (\ppp_tw + (H\cdot\nabla w))Bw dxdt.
                                                  \eqno{(3.2)}
$$
Now, we apply integration by parts to have:
\begin{align*}
& -2s \int_Q (\ppp_tw + (H\cdot \nabla w))Bw dxdt \\
=& -2s\int_Q (\ppp_tw)Bw dxdt - 2s\int_Q (H\cdot \nabla w) Bw dxdt
=: J_1 + J_2.
\end{align*}
Then 
$$
J_1 = -s\int_Q B\ppp_t(w^2) dxdt
= -s\int_{\OOO} \left[ B(x,t)\vert w(x,t)\vert^2\right]^T_0 dx
+ s\int_Q (\ppp_tB) w^2 dxdt.
$$
Moreover,
\begin{align*}
& J_2 = \sum_{i=1}^N -2s\int_{Q_i} \sum_{k=1}^N 
h_k(\ppp_kw_i)B_iw_i dxdt 
= -s\sum_{i=1}^N \int_{Q_i} \sum_{k=1}^N 
h_k B_i\ppp_k(\vert w_i\vert^2) dxdt \\
= -& s\sum_{i=1}^N \int_{\ppp\OOO_i \times (0,T)} 
(H\cdot \nu_{\OOO_i}) B_i\vert w_i\vert^2 dS_i dt 
+ s\sum_{i=1}^N \int_{Q_i} \mbox{div}\, (B_iH) \vert w_i\vert^2 dxdt.
\end{align*}
Here and henceforth $\nu_{\OOO_i}$ and $\nu_{\OOO}$ denotes the 
unit outward normal vectors to $\ppp\OOO_i$ and $\ppp\OOO$ respectively.
Then, recalling that $\nu_{i\to j}$ denotes the unit normal vector to 
$\gamma_{ij}:= \ppp\OOO_i \cap \ppp\OOO_j$ directed from $\OOO_i$ to $\OOO_j$,
we note that $\nu_{\OOO_i} = \nu_{i \to j}$ if $\OOO_i$ touches
$\OOO_j$ with $j\ne i$.
Moreover, let $dS_i$ denote the surface element of $\ppp\OOO_i$.

Therefore, (3.1) implies
\begin{align*}
& \sum_{i=1}^N s^2\int_{Q_i} B_i(x,t)^2 \vert w_i\vert^2 dxdt 
+ s\sum_{i=1}^N \int_{Q_i} (\mbox{div}\, (B_iH) - \vert\beta\vert)
 \vert w_i\vert^2 dxdt \\
+& s\sum_{i=1}^N \int_{\OOO_i} B_i(x,0) \vert w_i(x,0)\vert^2 dx
-  s\sum_{i=1}^N \int_{\ppp\OOO_i \times (0,T)} (H\cdot \nu_{\OOO_i})
B_i \vert w_i\vert^2 dS_idt 
\end{align*}
$$
\le C\sum_{i=1}^N \int_{Q_i} \vert P_iw_i\vert^2 dxdt 
+ s\sum_{i=1}^N \int_{\OOO_i} \vert B_i(x,T)\vert \vert w_i(x,T)\vert^2 dx.
                                          \eqno{(3.3)}
$$

Now we choose $\va_i$, that is, constants $r_i>0$ and $\beta > 0$ in the 
weight function $\va_i(x,t)$.  We set 
$R:= \max_{x\in \ooo{\OOO}} \vert x\vert$.
Thanks to (ii) of Condition (B), we can find a constant $\delta>0$ such that 
$$
\min_{1\le i \le N}\min_{x\in \ooo{\OOO_i}} (v_i\cdot H_i(x)) 
=: \delta > 0.
$$
Moreover, 
\begin{align*}
& B_i(x,t) = (\nabla \va_i(x,t) \cdot H(x)) - \beta
= (2(x+r_iv_i)\, \cdot\, H(x)) - \beta\\
=& 2r_i(v_i\cdot H(x)) + 2(x\cdot H(x)) - \beta 
\ge 2r_i\delta - 2R\Vert H\Vert_{C(\ooo{\OOO})} - \beta.
\end{align*}
We choose $r_i > 0$ for $1\le i \le N$ sufficiently large and
$\beta > 0$ sufficiently small such that 
$$
\min_{1\le i \le N} r_i 
> \frac{R\Vert H\Vert_{C(\ooo{\OOO})} }{\delta} + 1
                                               \eqno{(3.4)}
$$
and
$$
0 < \beta < \delta.                             \eqno{(3.5)}
$$
Then, noting the definition of $B$ in (3.1), we have
$$
\min_{1\le i \le N}\min_{(x,t)\in \ooo{Q_i}} B_i(x,t) 
= \min_{(x,t)\in \ooo{Q}} B(x,t) > \delta.      \eqno{(3.6)}
$$
Henceforth $C>0$ denotes generic constants dependent on 
$v_i, r_i, 1\le i \le N$ and $H, \beta, \OOO, T$, but independent of 
$s>0$.  Hence the application of (3.6) in (3.3) yields 
\begin{align*}
& s^2\sum_{i=1}^N \int_{Q_i} \vert w_i\vert^2 dxdt 
- Cs\sum_{i=1}^N \int_{Q_i} \vert w_i\vert^2 dxdt 
+ s\sum_{i=1}^N \int_{\OOO_i} \vert w_i(x,0)\vert^2 dx\\
- & s\sum_{i=1}^N \int_{\ppp\OOO_i \times (0,T)} (H\cdot \nu_{\OOO_i})
B_i \vert w_i\vert^2 dS_idt \\
\le & C \int_Q \vert Pw\vert^2 dxdt 
+ Cs \sum_{i=1}^N \int_{\OOO_i} \vert w_i(x,T)\vert^2 dx.
\end{align*}
Choosing $s_0 > 0$ sufficiently large, we can absorb the second term on the 
left-hand side into the first term, we obtain
$$
s^2\int_Q \vert w\vert^2 dxdt + s\int_{\OOO} \vert w(x,0)\vert^2 dx
+ \left(- s\sum_{i=1}^N \int_{\ppp\OOO_i \times (0,T)} (H\cdot \nu_{\OOO_i})
B_i \vert w_i\vert^2 dS_idt \right)
$$
$$
\le C \int_Q \vert Pw\vert^2 dxdt 
+ Cs \int_{\OOO} \vert w(x,T)\vert^2 dx                          \eqno{(3.7)}
$$
for all large $s \ge s_0$.
\\
{\bf Second Step: estimation of the integral on 
$(\ppp\OOO_i \cap \OOO) \times (0,T)$}.
\\
We have
$$
-s\sum_{i=1}^N \int_{\ppp\OOO_i \times (0,T)} (H\cdot \nu_{\OOO_i})
B_i \vert w_i\vert^2 dSdt 
$$
$$
= -s\sum_{i=1}^N \int_{(\ppp\OOO_i \cap \OOO)\times (0,T)} 
(H\cdot \nu_{\OOO_i}) B_i \vert w_i\vert^2 dSdt 
- s\sum_{i=1}^N \int_{(\ppp\OOO_i\cap \ppp\OOO) \times (0,T)} 
(H\cdot \nu_{\OOO_i})B_i \vert w_i\vert^2 
dSdt.                                                   \eqno{(3.8)}
$$
Then, the integrals on the interfaces:
$$
-s\sum_{i=1}^N \int_{(\ppp\OOO_i \cap \OOO)\times (0,T)} 
(H\cdot \nu_{\OOO_i}) B_i \vert w_i\vert^2 dSdt 
$$
should be non-negative.  To this end, we have to choose $\va_i(x,t)$ and 
a partition $\OOO_i$ for $1 \le i \le N$.  It suffices to consider only 
$\ppp\OOO_i \cap \OOO$ whose measure is positive. 
Then, an interface $\ppp\OOO_i \cap \OOO$ is divided as
$$
\ppp\OOO_i \cap \OOO = \bigcup_j \gamma_{ij},
$$
where $j$ varies over a set which satisfies $\vert \gamma_{ij}\vert \ne 0$ and
$(H\cdot \nu_{i\to j}) \ne 0$ on $\gamma_{ij}$ for each 
$i \in \{1, ..., N\}$.
Then,
\begin{align*}
& - s\sum_{i=1}^N \int_{(\ppp\OOO_i\cap \OOO) \times (0,T)} 
(H\cdot \nu_{\OOO_i})B_i \vert w_i\vert^2 dSdt\\
= & -s\sum_{i=1}^N \sum_{\vert \gamma_{ij}\vert \ne 0, 
(H\cdot \nu_{i\to j}) \ne 0} \int_{\gamma_{ij}\times (0,T)}
  (H\cdot \nu_{i \to j})B_i \vert w_i\vert^2 dSdt.
\end{align*}
Using $(H\cdot \nu_{i\to j}) = -(H\cdot \nu_{j\to i})$, in the above
summation on the right-hand side, we can pair the terms
$$
\int_{\gamma_{ij}\times (0,T)} (H\cdot \nu_{i \to j})B_i \vert w_i\vert^2 dSdt
\quad \mbox{and}
\int_{\gamma_{ij}\times (0,T)} (H\cdot \nu_{j \to i})B_j \vert w_j\vert^2 dSdt,
$$
so that 
we can represent $- s\sum_{i=1}^N \int_{(\ppp\OOO_i\cap \OOO) \times (0,T)}
(H\cdot \nu_{\OOO_i})B_i \vert w_i\vert^2 dSdt$ by the sum of the terms 
$S_{ij}$ where $(i,j)$ varies over some subset $\Psi$ of
$\{ (i,j)\in \{1,..., N\}^2;\, i<j\}$:
\begin{align*}
& S_{ij}:= -s\left( \int_{\gamma_{ij}\times (0,T)} 
(H\cdot \nu_{i\to j}) B_i \vert w_i\vert^2 dSdt
+  \int_{\gamma_{ij}\times (0,T)} 
(H\cdot \nu_{j\to i}) B_j \vert w_j\vert^2 dSdt \right) \\
= & s\int_{\gamma_{ij}\times (0,T)} 
(H\cdot \nu_{j\to i}) (B_i \vert w_i\vert^2 - B_j \vert w_j\vert^2) dSdt.
\end{align*}
Since $u\in H^1(Q)$, we have $w_i=ue^{s\va_i}$ and $w_j = ue^{s\va_j}$ on 
$\gamma_{ij}$ by the trace theorem to $u$ on $\gamma_{ij}$,
we obtain 
$$
 - s\sum_{i=1}^N \int_{(\ppp\OOO_i\cap \OOO) \times (0,T)} 
(H\cdot \nu_{\OOO_i})B_i \vert w_i\vert^2 dSdt
= \sum_{(i,j) \in \Psi}S_{ij} 
$$
$$
=  s\sum_{(i,j) \in \Psi} 
\int_{\gamma_{ij}\times (0,T)} (H\cdot \nu_{j \to i}) 
(B_ie^{2s\va_i} - B_je^{2s\va_j})\vert u\vert^2 dSdt.  \eqno{(3.9)}
$$
In order to obtain the non-negativity of the sum of the integrals in (3.9), 
we prove    
\\
{\bf Lemma 3.1.}
\\
{\it
Under Condition (C), we can choose $r_1, ..., r_N > 0$ such that 
$$
r_i^2 > 4r_j^2 + 6R^2 \quad \mbox{for $(i,j) \in \Psi$}     \eqno{(3.10)}
$$
and 
$$
\min_{1\le i\le N} r_i > \max\left\{ \frac{R\Vert H\Vert_{C(\ooo{\OOO})}}
{\delta} + 1, \, R\right\}.                                    \eqno{(3.11)}
$$
}

Postponing the proof of Lemma 3.1, we continue the proof of Theorem 2.1.
By means of (3.10), we can find a constant $\delta_2 > 0$ such that 
$$
r_i^2 - 4r_j^2 - 6R^2 > \delta_2 \quad \mbox{for all $(i,j) \in \Psi$}.
                                                       \eqno{(3.12)}
$$
We set $r^*:= \max_{1\le i \le N} r_i$.  

We choose $s_1 = s_1(\delta_2,r^*, \OOO, T, H, \beta) > 0$ such that
$s_1 > s_0$ and 
$$
e^{s_1\delta_2} > \frac{2(r^*+R)\Vert H\Vert_{C(\ooo\OOO)} + \beta}
{\delta}.                                             \eqno{(3.13)}
$$
By (3.1) and the fact that $H \in C(\ooo{\OOO})$, we see that 
$$
\vert B_j(x,t)\vert
= \vert 2r_j(v_j\cdot H(x)) + 2(x\cdot H(x)) - \beta \vert 
\le (2r_j+2R)\Vert H\Vert_{C(\ooo{\OOO})} + \beta
$$
for $x \in \ooo{\OOO}$ and $0<t<T$.
Therefore,
$$
 B_{i}e^{2s\va_i} - B_je^{2s\va_j}
\ge \delta e^{2s\va_{i}} - (2(r_{j}+R)\Vert H\Vert_{C(\ooo{\OOO})}
+ \beta) e^{2s\va_j}
$$
$$
= e^{2s\va_i} (\delta - (2(r_j+R)\Vert H\Vert_{C(\ooo{\OOO})}
+ \beta) e^{-2s(\va_i-\va_j)}.
                                          \eqno{(3.14)}
$$
On the other hand, using the inequality
$r_i - \vert x\vert \ge r_i - R > 0$, by (3.11) and the fact that  
$$
\vert x+r_jv_j\vert^2 \le 2\vert x\vert^2 + 2\vert r_jv_j\vert^2,
\quad (r_i - R)^2 \ge \frac{1}{2}r_i^2 - R^2,
$$
we have that 
\begin{align*}
& \va_i(x,t) - \va_j(x,t) 
= \vert x+v_ir_i\vert^2 - \vert x+v_jr_j\vert^2\\
\ge& (\vert v_ir_i\vert - \vert x\vert)^2
- 2(\vert x\vert^2 + \vert v_jr_j\vert^2)
\ge (r_i - R)^2 - 2(R^2 + r_j^2)\\
\ge& \frac{1}{2}r_i^2 - R^2 - 2(R^2 + r_j^2)
= \frac{1}{2}r_i^2 - 2r_j^2 - 3R^2 \ge \delta_2,
\end{align*}
where $\delta_2 > 0$ is independent of $(i, j) \in \Psi$.
Therefore, 
$$
 -2s(\va_i(x,t) - \va_j(x,t))
\le -s(r_i^2 - 4r_j^2 - 6R^2) \le -s\delta_2,
$$
and so
$$
 e^{-2s(\va_i(x,t) - \va_j(x,t))} \le e^{-s\delta_2}.
$$
Hence, (3.13) implies 
$$
\delta - (2(r_j+R)\Vert H\Vert_{C(\ooo\OOO)} + \beta)
e^{-2s(\va_i(x,t) - \va_j(x,t))}
$$
$$
\ge \delta - e^{-s_1\delta_2}(2(r^*+R)\Vert H\Vert_{C(\ooo\OOO)} + \beta)
> 0 \quad \mbox{for all $s \ge s_1$}.               \eqno{(3.15)}
$$
Consequently (3.14) and (3.15) yield
\begin{align*}
& (B_ie^{2s\va_i} - B_je^{2s\va_j})(x,t)
\ge e^{2s(\vert x_i + r_iv_i\vert^2 - \beta t)}
(\delta - e^{-s_1\delta_2}(2(r^*+R)\Vert H\Vert_{C(\ooo\OOO)} + \beta))\\
\ge& e^{-2s\beta T}(\delta - e^{-s_1\delta_2}
(2(r^*+R)\Vert H\Vert_{C(\ooo\OOO)} + \beta))
\ge 0 \quad \mbox{for all $s \ge s_1$}.
\end{align*}
Thus, we obtain 
$$
s\sum_{(i,j) \in \Psi} \int_{\gamma_{ij}\times (0,T)}
(H\cdot \nu_{j\to i})
(B_ie^{2s\va_i} - B_je^{2s\va_j})(x,t)\vert u\vert^2 dSdt \ge 0
\quad \mbox{for all $s\ge s_1$}.                 \eqno{(3.16)}
$$
Moreover, we have that 
$$
 -s\sum_{i=1}^N \int_{(\ppp\OOO_i\cap \ppp\OOO)\times (0,T)} 
(H\cdot \nu_{\OOO}) B\vert w\vert^2 dS dt
$$
$$
= s\sum_{i=1}^N \int_{\ppp\OOO_- \times (0,T)} (H\cdot \nu_{\OOO}) B
\vert w\vert^2 dS dt
- s\int_{\ppp\OOO_+ \times (0,T)} (H\cdot \nu_{\OOO}) B\vert w\vert^2 dS dt.
                                          \eqno{(3.17)}
$$
Therefore, (3.8), (3.9), (3.16) and (3.17) imply 
$$
 -s\sum_{i=1}^N \int_{\ppp\OOO_i \times (0,T)}
(H\cdot \nu_{\OOO_i})B_i\vert w_i\vert^2 dS_idt
$$
$$
\ge s\int_{\ppp\OOO_- \times (0,T)} (H\cdot \nu_{\OOO}) B\vert w\vert^2 dS dt
- s\int_{\ppp\OOO_+ \times (0,T)} (H\cdot \nu_{\OOO}) B\vert w\vert^2 dS dt.
                                              \eqno{(3.18)}
$$
Substituting (3.18) into (3.7), we can complete the proof of Theorem 2.1.
\\
\vspace{0.2cm}
\\
{\bf Third Step: proof of Lemma 3.1.}
\\
Before the proof, we provide a heuristic explanation. 
In terms of directed graph, we can interpret the condition  
$(H\cdot \nu_{j\to i}) > 0$ saying that $O_i$ is a downstream node
of $O_j$ with respect to the stream field $H$.  
Let the directed graph $\Lambda$ be composed of nodes $O_1, ..., O_N$ 
and directed segments 
$\Gamma_{ij}$ connecting $O_i$ and $O_j$ where $i$ and $j$ vary over some 
subset of $\{ 1, ..., N\}^2$. Then, we assign $r_j$ to $O_j$ for each 
$j \in \{ 1, ..., N\}$.
Condition (3.10) means that we have to choose constants $r_i > 0$
sufficiently large 
compared with the values $r_j$ assigned to all the upstream nodes $O_j$
in view of $O_i$, and such assignments of $r_i$ must be made over all the nodes
connected mutually.
We note that we can assign $r_i$ independently of $r_j$ if 
$(H\cdot \nu_{j\to i}) = 0$ on $\gamma_{ij}$.
To this end, the ordering of the assignments of the value $r_i$ to $O_i$ 
is essential.
We can see that we should choose $r_i$ at a downstream node
after the assignments of $r_j$ at all the upstream nodes.
The assignment procedure cannot be limited only to related nodes but 
is concerned with the structure of the directed graph $\Lambda$.
We prove Lemma 3.1 by induction concerning the number of the nodes of
$\Lambda$.
\\
Now we can complete
\\
{\bf Proof of Lemma 3.1.}
\\
First we introduce some notations.
Let $\Lambda$ be a directed graph composed of nodes
$O_1, ..., O_N$ and the directed segments connecting some of 
the nodes.  We can assume that $\Lambda $ has a directed segment, that is,
there exists at least one node connected to other nodes.  Otherwise we can 
choose $r_i$ for $1\le i \le N$ satisfying only (3.10),
and so the proof of Lemma 3.1 is trivial.

For $O_i$ for $1 \le i \le N$, we denote the directed segments
at $O_i$ as follows: We choose disjoint sets $J_+(O_i), J_-(O_i) \subset 
\{ 1, ..., N\}$ such that all the directed segments at $O_i$ 
are denoted by $\Gamma_{ij}$ with $j \in J_+(O_i) \cup J_-(O_i)$, and 
$(H\cdot \nu_{i\to j}) > 0$ on $\gamma_{ij}$ if $j\in J_+(O_i)$
and $(H\cdot \nu_{i\to j}) < 0$ on $\gamma_{ij}$ if $j\in J_-(O_i)$.
We note that $J_+(O_i)$ or $J_-(O_i)$ may be an empty set.
Then,
\\
{\bf Definition 3.1.}
\\
{\it
We call a node $O_i$ a terminus node of the directed graph $\Lambda$
if $J_+(O_i) = \emptyset$ and $J_-(O_i) \ne \emptyset$.
}

The definition means that 
there are no outgoing segments from $O_i$.
Then,
\\
{\bf Lemma 3.2.}
\\
{\it
(i) If $\Lambda$ has no closed loops, then there are no closed loops
in any directed graph $\www{\Lambda}$ obtained by deleting $O_i$ and 
$\Gamma_{ij}$ with $j\in J_+(O_i) \cup J_-(O_i)$ from $\Lambda$.
\\
(ii) If $\Lambda$ has no closed loops, then $\Lambda$ has a terminus node.
}
\\
{\bf Proof of Lemma 3.2.}
\\
{\bf (i)}.
As is seen by the construction $\www{\Lambda}$, the directed segments of
$\www{\Lambda}$ is a subset of the ones of $\Lambda$.  Since the directed 
segments do not form a closed loop, the subset of the totality of the 
directed segments cannot create closed loops.
Thus the proof of Lemma 3.2 (ii) is complete.
\\
{\bf (ii)}.  
First, 
by deleting the nodes of $\lambda$ which are not connected to any other nodes,
we construct $\Lambda_1 \subset \Lambda$ satisfying 
$J_+(O_i) \cup J_-(O_i) \ne \emptyset$ for each node $O_i$ of $\Lambda_1$.

If $\Lambda_1$ has a terminus node, then the proof is already finished.
Hence, assume that $\Lambda_1$ has no terminus nodes.  Then,  
$$
\mbox{for each $i \in \{ 1, ..., N\}$, 
there exists an outgoing segment from $O_i$.}                    \eqno{(3.19)}
$$
We start at an arbitrarily chosen node of $\Lambda_1$ and move to other nodes
along arbitrarily chosen directed outgoing segments from the start node.
In view of (3.19), we can find an outgoing segment as an exit at any node 
where we enter.  Therefore, the tour can be continued endlessly.
In particular, during this node-to-node $(N+1)$-times movements, 
we have to pass the same node at least twice, because 
the total number of the nodes is $N$ and the movements can be 
continued an infinite times.  This means that a path with 
any $(N+1)$-movement contains a closed loop.  
Thus the proof of (ii) is complete as well as that Lemma 3.2.
$\blacksquare$

Now we proceed to
\\
{\bf Completion of the proof of Lemma 3.1.}
Our proof is by the induction with respect to the numbers $N \ge 2$ 
of nodes of directed graphs.
\\
We first prove Lemma 3.1 for any directed graph with $N=2$.
We can have two cases\\
(i) $(H\cdot \nu_{2 \to 1}) = 0$ on $\gamma_{12}$.
\\
(ii) $(H\cdot \nu_{2\to 1}) > 0$ on $\gamma_{12}$.

The case $(H\cdot \nu_{2\to 1}) < 0$ on $\gamma_{12}$ is reduced to 
(ii) by exchanging $O_1$ and $O_2$.

In case (i), we do not need to consider (3.10) and we can satisfy 
(3.11) by choosing sufficiently large $r_1, r_2$.
In case (ii), (3.10) is reduced to $r_1^2 > 4r_2^2 + 6R^2$.
Therefore, we can immediatetly see that there exist $r_1, r_2 > 0$
satisfying (3.10) and (3.11).  Thus Lemma 3.1 is proved for $N=2$.

We assume that Lemma 3.1 is proved for arbitrary directed graphs $\Lambda$
with an $N$-number of nodes.
We will prove Lemma 3.1 for arbitrary directed graphs $\www{\Lambda}$ with 
$(N+1)$-number nodes which have no closed loops.
By Lemma 3.2 (ii), it follows that $\www{\Lambda}$ has a terminus node, which 
is denoted by $O_{N+1}$.  Then we can find a non-empty set 
$J_-(O_{N+1}) \subset \{ 1, ..., N\}$ such that 
each directed segment incoming to $O_{N+1}$ is labelled by 
$\Gamma_{N+1,j}$ with $j\in J_-(O_{N+1})$ and 
$(H\cdot \nu_{N+1\to j}) < 0$ on $\gamma_{N+1,j}$ for $j\in J_-(O_{N+1})$.
We construct another directed graph $\widehat{\Lambda}$ by deleting $O_{N+1}$
and $\Gamma_{N+1,j}$ for all $j\in J_-(O_{N+1})$.
Then, $\widehat{\Lambda}$ has $N$-nodes and by Lemma 3.1 (i), we see that 
$\widehat{\Lambda}$ has no closed loops.
Since the conclusion holds for $\widehat{\Lambda}$ with $N$-nodes 
by the assumption of the induction, 
Lemma 3.1 is satisfied 
for $\widehat{\Lambda}$: there exist constants $r_1, ..., r_N > 0$ such that 
(3.10) and (3.11) hold.  Since $r_1, ..., r_N$ are already determined, it 
suffices to choose $r_{N+1} > 0$ such that
$$
r_{N+1}^2 >\max_{1\le j \le N} (4r_j^2 + 6R^2)
$$
and 
$$
r_{N+1} > \max\left\{ \frac{R\Vert H\Vert_{C(\ooo\OOO)}}{\delta} + 1,\,\,
R \right\}.
$$
Then, we can readily verify that (3.10) and (3.11) are satisfied for 
$\Lambda$ with $(N+1)$-nodes.  Therefore, we proved that Lemma 3.1 holds 
true also for arbitrary $(N+1)$-node directed graphs without closed 
loops.  Thus the induction finishes the proor of Lemma 3.1, which 
completes the proof of Theorem 2.1.
$\blacksquare$
\section{Applications to inverse problems}

{\bf 4.1. Observability inequality.}

In this section, we apply Theorem 2.1 to two kinds of inverse problems
for first-order transport equations.
Throughout this section, we assume Conditions (B) and (C) stated in 
Section 2.

We consider 
$$
\ppp_tu(x,t) + (H(x)\cdot \nabla u(x,t)) + p(x)u = 0
\quad \mbox{in $Q$}.                     \eqno{(4.1)}
$$
Then,
\\
{\bf Theorem 4.1 (observability inequality).}\\
{\it
We assume Conditions (B) and (C) stated in Section 2.
There exists a constant $T_0 > 0$ depending on $H, \OOO$, such that
for all $T > T_0$, we can find a constant $C>0$ such that 
$$
\Vert u(\cdot,0)\Vert_{L^2(\OOO)} \le 
C\Vert u\Vert_{L^2(\ppp\OOO
\times (0,T))}                         \eqno{(4.2)}
$$
for all $u\in H^1(Q)$ satisfying (4.1).
Here the constant $C>0$ depends on $\OOO, T, \va$.
}
\\

By Theorem 2.1, we can prove other types of stability estimates of solution 
$u$ to for (4.1) such as an estimate of H\"older type 
(e.g., Chapter 3 in Yamamoto \cite{Y24}), but we here 
omit further details.
The proof is derived from a relevant Carleman estimate (Theorem 2.1), and 
here we provide the proof in a more general setting. 
Assume that we know a Carleman 
estimate: there exist a constant $\beta > 0$ and a function $d(x)$ 
which is of piecewise $C^1$ on $\ooo{\OOO}$, such that 
$$
s\int_{\OOO} \vert u(x,0)\vert^2 e^{2s\va(x,0)} dx 
+ s^2\int_Q \vert u\vert^2 e^{2s\va} dxdt 
$$
$$
\le C\int_Q \vert \ppp_tu + (H\cdot \nabla u) + pu\vert^2 e^{2s\va} dxdt 
+ C\int_{\ppp\OOO_+ \times (0,T)} \vert u\vert^2 e^{2s\va} dSdt
$$
$$
+ Cs\int_{\OOO} \vert u(x,T)\vert^2 e^{2s\va(x,T)} dx 
\quad \mbox{for all large $s>0$}.                       \eqno{(4.3)}
$$
Here we set $\va(x,t):= d(x)- \beta t$ for $(x,t) \in Q:= \OOO\times 
(0,T)$ and recall that 
$\ppp\OOO_+:= \{ x\in \ppp\OOO;\, (H(x)\cdot \nu(x)) \ge 0\}$.

For the proof of Theorem 4.1, it suffices to prove
\\
{\bf Proposition 4.1.}
\\
{\it 
We assume Carleman estimate (4.3).  Then the conclusion of Theorem 4.1 
holds true.}
\\
{\bf Proof of Proposition 4.1.}
\\
{\bf First Step.}
\\
We can prove an energy estimate.\\
{\bf Lemma 4.1.}\\
{\it
Let $u\in H^1(Q)$ satisfy 
$$
\ppp_tu(x,t) + (H(x)\cdot \nabla u(x,t)) + p(x)u(x,t) = F(x,t), \quad 
(x,t) \in Q                            \eqno{(4.4)}
$$
with $\Vert p\Vert_{L^{\infty}(Q)} \le M$.
Then there exists a constant $C>0$ which is dependent only on
$M$, $\OOO$ and $T$, $H$, such that 
\begin{align*}
& \int_{\OOO} \vert u(x,t)\vert^2 dx
+ \Vert u\Vert^2_{L^2(\ppp\OOO_+ \times (0,t))}\\
\le & C(\Vert F\Vert^2_{\LLLL} 
+ \Vert u\Vert^2_{L^2(\ppp\OOO_-\times (0,T))}
+ \Vert u(\cdot,0)\Vert^2_{L^2(\OOO)}), \quad 0\le t\le T.
\end{align*}
}
\\
{\bf Proof of Lemma 4.1.}\\
We recall that $H = (h_1, ..., h_d)$.
Multiplying (4.4) with $2u(x,t)$ and integrating over $\OOO$, we have
$$
\int_{\OOO} 2u(x,t)\ppp_tu(x,t)dx 
+ \int_{\OOO} 2u(x,t)\sum^d_{j=1} h_j\ppp_ju dx
+ \int_{\OOO} 2p\vert u(x,t)\vert^2 dx = \int_{\OOO} 2F(x,t)u(x,t)dx,
$$
that is,
$$
\ppp_t\int_{\OOO} \vert u(x,t)\vert^2 dx 
+ \int_{\OOO} (H\cdot \nabla (u^2)) dx
+ \int_{\OOO} 2p u^2 dx = \int_{\OOO} 2Fu \,dx
$$
in $(0,T)$.
Setting $E(t):= \int_{\OOO} \vert u(x,t)\vert^2 dx$ and integrating by 
parts, by $\Vert p\Vert_{L^{\infty}(Q)} \le M$ we obtain
\begin{align*}
& E'(t) - \int_{\OOO} (\ddd H)\vert u(x,t)\vert^2 dx
+ \int_{\ppp\OOO_+} (H\cdot\nu) u^2 dS
+ \int_{\ppp\OOO_-} (H\cdot\nu) u^2 dS\\
=& \int_{\OOO} 2Fu\, dx - 2\int_{\OOO} p u^2 dx
\le CME(t) + \int_{\OOO} (\vert F\vert^2 + u^2) dx.
\end{align*}
Therefore, using $H \in C^1(\ooo{\OOO})$ and 
rewriting in $\eta$, we have
\begin{align*}
& E'(\eta) + \int_{\ppp\OOO_+} \vert u(x,\eta)\vert^2 dS\\
\le& CE(\eta) 
+ C\int_{\ppp\OOO_-} \vert u(x,\eta)\vert^2 dS
+ C\int_{\OOO} \vert F(x,\eta)\vert^2 dx, \quad 0\le \eta \le T.
\end{align*}
Hence, integrating in $\eta \in (0,t)$, we obtain
\begin{align*}
& E(t) + \int^t_0 \left( \int_{\ppp\OOO_+} \vert u(x,\eta)\vert^2 dS\right)
d\eta\\
\le& E(0) + C\int^t_0 E(\eta) d\eta 
+ C(\Vert u\Vert^2_{L^2(\ppp\OOO_- \times (0,T))}
+ \Vert F\Vert_{L^2(Q)}^2), \quad 0\le t\le T.
\end{align*}
The Gronwall inequality completes the proof of Lemma 4.1.
$\blacksquare$
\\
{\bf Second Step.}
\\
We choose $T_0 > 0$ such that 
$$
T_0 > \frac{1}{\beta}(\max_{x\in \ooo{\OOO}} d(x)
- \min_{x\in \ooo{\OOO}} d(x))                      \eqno{(4.5)}
$$
and we assume that $T>T_0$.
We set 
$$
\mu:= \min_{x\in \ooo{\OOO}} d(x) + \beta T - \max_{x\in \ooo{\OOO}} d(x) > 0
                                                 \eqno{(4.6)}
$$
for $T > T_0$.
Applying (4.3) to (4.1), we have
$$
 \int_{\OOO} \vert u(x,0)\vert^2 e^{2s\va(x,0)} dx 
\le C\int_{\ppp\OOO_+ \times (0,T)} \vert u\vert^2 e^{2s\va} dSdt
+ C\int_{\OOO} \vert u(x,T)\vert^2 e^{2s\va(x,T)} dx
$$
for all large $s>0$.

Therefore, applying Lemma 4.1, we obtain
\begin{align*}
& e^{2s\min_{x\in \ooo{\OOO}} \va(x,0)}\Vert u(\cdot,0)\Vert_{L^2(\OOO)}^2\\
\le& Ce^{Cs}\Vert u\Vert^2_{L^2(\ppp\OOO_+ \times (0,T))}
+ Ce^{2s\max_{x\in \ooo{\OOO}} \va(x,T)}\Vert u(\cdot,T)\Vert^2_{L^2(\OOO)}\\
\le& Ce^{Cs}\Vert u\Vert^2_{L^2(\ppp\OOO_+ \times (0,T))}
+ Ce^{2s\max_{x\in \ooo{\OOO}} \va(x,T)} \Vert u\Vert^2
_{L^2(\ppp\OOO_- \times (0,T))}
+ Ce^{2s\max_{x\in \ooo{\OOO}} \va(x,T)}\Vert u(\cdot,0)\Vert^2
_{L^2(\OOO)}
\end{align*}
for all large $s>0$.  Dividing by 
$e^{2s\min_{x\in \ooo{\OOO}} \va(x,0)}$ we reach 
$$
\Vert u(\cdot,0)\Vert_{L^2(\OOO)}^2
\le Ce^{Cs}\Vert u\Vert^2_{L^2(\ppp\OOO \times (0,T))}
+ Ce^{-2s\mu}\Vert u(\cdot,0)\Vert^2_{L^2(\OOO)}
$$
for all large $s>0$.
Choosing $s>0$ sufficiently large, we can absorb the second term 
on the right-hand side into the left-hand side thanks to (4.6).
Thus, the proof of Proposition 4.1, and the one of Theorem 4.1 are complete.
$\blacksquare$
\\
{\bf 4.2. Inverse source problem}

We consider 
$$
\left\{ \begin{array}{rl}
& \ppp_tu(x,t) + (H(x)\cdot \nabla u(x,t)) + p(x)u(x,t) = R(x,t)f(x), \\
& u(x,0) = 0, \qquad x\in \OOO, \, 0<t<T.
\end{array}\right.                       \eqno{(4.7)}
$$
Let $p, R$ be given and $\Vert p\Vert_{L^{\infty}(Q)} \le M$, where 
$M>0$ is an arbitrarily fixed constant.

We discuss
\\
{\bf Inverse source problem\index{inverse source problem}.}\\
{\it 
Determine $f(x)$ for $x\in \OOO$ by data $u$ on 
$\ppp\OOO \times (0,T)$.}

By Theorem 2.1, we will prove
\\
{\bf Theorem 4.2.}\\
{\it
We assume Conditions (B) and (C) stated in Section 2, and  
$$
R, \, \ppp_tR \in L^2(0,T;L^{\infty}(\OOO)), \quad 
\vert R(x,0)\vert > 0, \quad x\in \ooo{\OOO}         \eqno{(4.8)}
$$
and $u, \ppp_tu \in H^1(Q)$.  Then, there exists $T_0 > 0$ depending on 
$H, \OOO$ such that for each $T > T_0$ we can find a constant $C>0$ 
such that 
$$
\Vert f\Vert_{L^2(\OOO)} \le C\Vert \ppp_tu\Vert
_{L^2(\ppp\OOO\times (0,T))}.                         \eqno{(4.9)}
$$
Here, the constant $C>0$ depends on $\OOO, T, H, M, \va, R$.
}

Similarly to Theorem 4.1, we can prove other types of stability estimates
by Theorem 2.1, but we concentrate on the global Lipschitz stability 
(4.9).

For the proof of Theorem 4.2, it suffices to prove
\\
{\bf Proposition 4.2.}
\\
{\it
We assume Conditions (B) and (C) stated in Section 2.
Let Carleman estimate (4.3) hold.  
Then, the conclusion of Theorem 4.2 holds
true.
}
\\
{\bf Proof of Proposition 4.2.}
\\
By (4.8), we see that $R \in H^1(0,T;L^{\infty}(\OOO)) \subset 
C([0,T];L^{\infty}(\OOO))$.  

We take the time-derivative of $u$: $y:= \ppp_tu$. Then
$$
\left\{ \begin{array}{rl}
& \ppp_ty(x,t) + (H(x)\cdot \nabla y(x,t)) + p(x)y(x,t) = (\ppp_tR)(x,t)f(x), \\
& y(x,0) = R(x,0)f(x), \qquad x\in \OOO, \, 0<t<T.
\end{array}\right.                       \eqno{(4.10)}
$$
In view of (4.8), the application of Lemma 4.1 to (4.10) yields 
$$
\Vert y(\cdot,T)\Vert_{L^2(\OOO)} 
\le C(\Vert (\ppp_tR)f\Vert_{\LLLL}
+ \Vert y\Vert_{L^2(\ppp\OOO_- \times (0,T))} 
+ \Vert R(\cdot, 0)f\Vert_{L^2(\OOO)})
$$
$$
\le C(\Vert f\Vert_{L^2(\OOO)} 
+ \Vert \ppp_tu\Vert_{L^2(\ppp\OOO_- \times (0,T))}).              \eqno{(4.11)}
$$
Next we apply (4.3) to (4.10):
$$
 s\int_{\OOO} \vert y(x,0)\vert^2 e^{2s\va(x,0)} dx     \eqno{(4.12)}
$$
$$
\le C\int_Q \vert (\ppp_tR)f\vert^2 \weight dxdt
+ Cs\int_{\OOO} \vert y(x,T)\vert^2 e^{2s\va(x,T)} dx
+ Cs\int_{\ppp\OOO_+\times (0,T)} y^2 \weight dSdt
$$
for all large $s>0$.

Since $\va(x,t) \le \va(x,0)$ for $(x,t) \in \ooo{Q}$, by (4.8)
we have
\begin{align*}
& \int_Q \vert \ppp_tR(x,t)\vert^2\vert f(x)\vert^2 e^{2s\va(x,t)} dxdt
\le \int_{\OOO} \left( \int^T_0 \vert \ppp_tR(x,t)\vert^2 dt\right)
\vert f(x)\vert^2 e^{2s\va(x,0)} dxdt\\
\le & \sup_{x\in \OOO} \left( \int^T_0 \vert \ppp_tR(x,t)\vert^2 dt\right)
\int_{\OOO} \vert f(x)\vert^2 e^{2s\va(x,0)} dx \\
\le & \Vert \ppp_tR\Vert^2_{L^2(0,T;L^{\infty}(\OOO))}
\int_{\OOO} \vert f(x)\vert^2 e^{2s\va(x,0)} dx.
\end{align*}
Consequently, (4.12) yields
\begin{align*}
& s\int_{\OOO} \vert y(x,0)\vert^2 e^{2s\va(x,0)} dx\\
\le& C\int_{\OOO} \vert f\vert^2 e^{2s\va(x,0)} dx
+ Cs\int_{\OOO} \vert y(x,T)\vert^2dx e^{2s\max_{x\in \ooo{\OOO}} \va(x,T)}
\end{align*}
$$
+ Cse^{Cs}\Vert \ppp_tu\Vert^2_{L^2(\ppp\OOO_+ \times (0,T))}
                                                       \eqno{(4.13)}
$$
for all large $s>0$.
We have
$$
\int_{\OOO} \vert f(x)\vert^2 e^{2s\va(x,0)} dx 
\le C\int_{\OOO} \vert y(x,0)\vert^2 e^{2s\va(x,0)} dx 
$$
by the last condition in (4.8), and the substitution of this 
inequality into (4.13) leads to  
\begin{align*}
& s\int_{\OOO} \vert f(x)\vert^2 e^{2s\va(x,0)} dx\\
\le& C\int_{\OOO} \vert f(x)\vert^2 e^{2s\va(x,0)} dx
+ Cs\left( \int_{\OOO} \vert y(x,T)\vert^2dx \right)
e^{2s\max_{x\in \ooo{\OOO}} \va(x,T)}
+ Cse^{Cs}\Vert \ppp_tu\Vert^2_{L^2(\ppp\OOO_+ \times (0,T))}
\end{align*}
for all large $s>0$.
Choosing $s>0$ large, we can absorb the first term on the right-hand side into
the left-hand side, and  we obtain
\begin{align*}
& s\int_{\OOO} \vert f(x)\vert^2 e^{2s\va(x,0)} dx\\
\le& Cs\int_{\OOO} \vert y(x,T)\vert^2dx 
e^{2s\max_{x\in \ooo{\OOO}} \va(x,T)}
+ Cse^{Cs}\Vert \ppp_tu\Vert^2_{L^2(\ppp\OOO_+ \times (0,T))}
\end{align*}
for all large $s>0$.
Dividing by $s e^{2s\min_{x\in \ooo{\OOO}} \va(x,0)}$ and noting 
(4.6), we reach
$$
\Vert f\Vert^2_{L^2(\OOO)} 
\le Ce^{Cs}\Vert \ppp_tu\Vert^2_{L^2(\ppp\OOO_+ \times (0,T))}
+ C\int_{\OOO} \vert y(x,T)\vert^2 dx e^{-2s\mu}                   \eqno{(4.14)}
$$
for all large $s>0$.

Applying (4.11) to the second term on the right-hand side of (4.14), we have
$$
\Vert f\Vert^2_{L^2(\OOO)}
\le Ce^{Cs}\Vert \ppp_tu\Vert^2_{L^2(\ppp\OOO_+ \times (0,T))}\\
+ C\Vert f\Vert^2_{L^2(\OOO)}e^{-2s\mu} 
+ C\Vert \ppp_tu\Vert^2_{L^2(\ppp\OOO_- \times (0,T))}e^{-2s\mu}
$$
for all large $s>0$.  Since $\mu>0$, we choose $s>0$ sufficiently large, 
so that we absorb the second term on the right-hand side into the left-hand 
side.  Therefore, the proof of Proposition 4.2 is complete.
Thus the proof of Theorem 4.2 is complete.
$\blacksquare$
\section{Concluding remarks}
{\bf 5.1.}
We have considered a first-order transport equation
$$
\ppp_tu(x,t) + (H(x)\cdot \nabla u(x,t)) + p(x)u(x,t) = F(x,t) 
\quad \mbox{in $Q$},                 \eqno{(5.1)}
$$
and established a Carleman estimate under a more general condition 
on $H$ than the one in the existing works Cannarsa, Floridia, G\"olgeleyen and 
Yamamoto \cite{CFGY}, Gaitan and Ouzzane \cite{GO}, G\"olgeleyen and 
Yamamoto \cite{GY}, Chapter 3 in Yamamoto \cite{Y24}.
The key for such an estimate (Theorem 2.1) is the weight 
function which is piecewise $C^1$ in $x\in \ooo{\OOO}$, and the 
condition of the weight requires that the associated directed graph 
by the stream field $H(x)$ has no closed loops.

Our Carleman inequality produces stability estimates for inverse problems,
the derivation of which follows standard arguments, once the Carleman 
estimate (Theorem 2.1) is established.
\\
{\bf 5.2.}
In this article, we assume that $H$ in (5.1) does not depend on $t$,
but we can similarly discuss the case where $H(x,t)$ depends also on 
$t$.  In a future work, we will study the details.
\section{Appendix: Proofs of Propositions 1.1 and 2.1}

We will prove a slightly more general Carleman estimate
(Proposition 6.1) than Proposition 1.1.
Choosing $d \in C^1(\ooo{\OOO})$ and a constant $\beta > 0$, we define
$$
\va(x,t) := d(x) - \beta t, \quad B(x,t) := (\nabla d(x)\cdot H(x)) - \beta
$$
for $(x,t) \in Q$.
We consider
$$
\ppp_tu(x,t) + (H(x)\cdot \nabla u(x,t)) + p(x)u(x,t) = F(x,t) 
\quad \mbox{in $Q$},                 \eqno{(6.1)}
$$
where $H \in C^1(\ooo{\OOO})$, $\vert H(x)\vert > 0$ for all $x \in 
\ooo{\OOO}$, $p\in L^{\infty}(\OOO)$ and $F \in L^2(Q)$.
Then,
\\
{\bf Proposition 6.1.}
\\
{\it
We assume
$$
\min_{(x,t)\in \ooo{Q}} B(x,t) =: \delta_3 > 0.      \eqno{(6.2)}
$$
Then, there exist constants $C>0$ and $s_0>0$, which are dependent on 
$\va, M, H, \OOO, T$ but independent of $s$, such that
\begin{align*}
& s^2\int_Q u^2 \weight dxdt 
+ s\int_{\OOO} B(x,0)\vert u(x,0)\vert^2 e^{2s\va(x,0)} dx
+ s\int_{\ppp\OOO_- \times (0,T)} \vert (H\cdot \nu)\vert Bu^2 dSdt\\
\le &C\int_Q \vert F\vert^2\weight dxdt
+ Cs\int_{\ppp\OOO_+\times (0,T)} \vert (H\cdot \nu)\vert B u^2 e^{2s\va} dSdt
                                                                       \\
+ & Cs\int_{\OOO} B(x,T)\vert u(x,T)\vert^2 e^{2s\va(x,T)} dx
\end{align*}
for all $s\ge s_0$ and all $u \in H^1(Q)$ satisfying (6.1).
}
\\
{\bf Proof of Proposition 6.1.}
\\
The proof is similar to Proposition 2.1 of Chapter 3 in \cite{Y24}.
We provide it for completeness.
As discussed at the beginning of the proof of Theorem 2.1 in 
Section 3, it is sufficient to prove the proposition with 
$p\equiv 0$.  Similarly to the beginning of the proof of Theorem 2.1,
setting $w:= e^{s\va}$, we have that 
\begin{align*}
&Pw(x,t) := e^{s\va(x,t)}(\ppp_t + (H\cdot \nabla))e^{-s\va(x,t)}w(x,t))\\
=& \ppp_tw(x,t) + (H(x)\cdot \nabla w(x,t)) - sB(x,t)w(x,t).
\end{align*}
Therefore, by the same way as (3.2), we obtain 
$$
 \int_Q \vert Pw\vert^2 dxdt 
\ge s^2\Vert Bw\Vert^2
- 2s\int_Q (\ppp_tw + (H\cdot \nabla w)) Bw dxdt.            \eqno{(6.3)}
$$
Now, integration by parts yields 
\begin{align*}
& -2s \int_Q (\ppp_tw + (H\cdot \nabla w))Bw dxdt \\
=& -2s\int_Q (\ppp_tw)Bw dxdt - 2s\int_Q (H\cdot \nabla w) Bw dxdt
=: I_1 + I_2.
\end{align*}
Then 
$$
I_1 = -s\int_Q B\ppp_t(w^2) dxdt
= -s\int_{\OOO} \left[ B(x,t)\vert w(x,t)\vert^2\right]^T_0 dx
+ s\int_Q (\ppp_tB) w^2 dxdt
$$
and
\begin{align*}
& I_2 = -2s\int_Q \sum^d_{j=1} h_j(\ppp_jw)Bw dxdt
= -s\int_Q \sum^d_{j=1} h_jB\ppp_j(w^2) dxdt\\
=& -s\int_{\ppp\OOO\times (0,T)} \sum^d_{j=1} h_jB\nu_j w^2 dSdt
+ s\int_Q \sum^d_{j=1} \ppp_j(h_jB) w^2 dxdt.
\end{align*}
Since $H \in C^1(\ooo{\OOO})$ and $B \in C^1(\ooo{Q})$, we obtain
\begin{align*}
& I_1 + I_2
\ge -Cs\int_Q w^2 dxdt \\
-& s\left( \int_{\OOO} \left[B w^2\right]^T_0 dx
+ \int_{\ppp\OOO\times (0,T)} (H\cdot \nu)B w^2 dSdt\right).
\end{align*}
Hence, (6.3) implies
\begin{align*}
& \int_Q \vert Pw\vert^2 dxdt
\ge \left( s^2\min_{(x,t)\in \ooo{Q}} B(x,t)^2 - Cs\right)
\int_Q w^2 dxdt \\
-& s\left( \int_{\OOO} \left[B w^2\right]^T_0 dx
+ \int_{\ppp\OOO\times (0,T)} (H\cdot \nu)B w^2 dSdt\right).
\end{align*}
Then,  by choosing $s>0$ sufficiently large, we obtain
\begin{align*}
& s^2\int_Q w^2 dxdt 
+ s\int_{\OOO} B(x,0)w(x,0)^2 dx 
   - s\int_{\OOO} B(x,T)w(x,T)^2 dx \\
+& s\int_{\ppp\OOO_-\times (0,T)} \vert (H\cdot \nu)\vert B w^2 dSdt
   - s\int_{\ppp\OOO_+\times (0,T)} \vert (H\cdot \nu)\vert B w^2 dSdt\\
\le & C\int_Q \vert Pw\vert^2 dxdt
\end{align*}
for all $s>s_0$.  Substituting $w = e^{s\va}u$, we complete the proof of
Proposition 6.1.
$\blacksquare$
\\
{\bf Proof of Proposition 1.1.}
\\
In terms of (1.6) and (1.7), setting $d(x) := \vert x+rv\vert^2$ and
$R:= \max_{x\in \ooo{\OOO}} \vert x\vert$, we have
\begin{align*}
& (\nabla d(x)\cdot \, H(x)) - \beta
= 2(x+rv\,\cdot\, H(x)) - \beta\\
=& 2r(v\cdot H(x)) + 2(x\cdot H(x)) - \beta
\ge 2r\delta_1 - 2R\Vert H\Vert_{C(\ooo{\OOO})} - \beta > 0
\end{align*}
for all $x \in \ooo{\OOO}$.  Therefore, (6.2) is satisfied, and so 
by Proposition 6.1, we obtain Proposition 1.1.
$\blacksquare$
\\
{\bf Proof of Proposition 2.1.}
We set $d(x):= \rho(x)$.  Since $H = \nabla \rho$ in $\OOO$ and
$\vert \nabla \rho\vert > 0$ on $\ooo{\OOO}$, we see that 
$$
(\nabla d(x)\, \cdot \,H(x)) - \beta
= (\nabla \rho(x)\, \cdot \,\nabla\rho(x)) - \beta > 0
$$
for $x \in \ooo{\OOO}$, provided that $\beta > 0$ is chosen sufficiently
small.  Thus Proposition 6.1 completes the proof of 
Proposition 2.1.

{\bf Acknowledgements.}
Piermarco Cannarsa was supported by the PRIN 2022 PNRR Project P20225SP98 
"Some mathematical approaches to climate change and its impacts" 
(funded by the European Community-Next Generation EU), CUP E53D2301791 0001, 
by the INdAM  (Istituto Nazionale di Alta Matematica) Group for Mathematical 
Analysis, Probability and Applications, and by the MUR Excellence Department 
Project awarded to the Department of Mathematics, University of Rome 
Tor Vergata, CUP E83C23000330006.
Masahiro Yamamoto was supported by 
Grant-in-Aid for Scientific Research (A) 20H00117 
and Grant-in-Aid for Challenging Research (Pioneering) 21K18142 of 
Japan Society for the Promotion of Science.
Most of this work has been done when Masahiro Yamamoto was 
a visiting professors at University of Rome Tor Vergata and
at Sapienza University of Rome in 2024 and 2023.

\end{document}